 \documentclass[10pt]{article}
 \usepackage[top=1in,bottom=1in,left=1.5in,right=1.5in]{geometry}
  \usepackage{amsmath,amssymb}
  \usepackage{latexsym}
  \usepackage[dvips]{pstcol} 
  \usepackage[dvips]{graphicx}

  \newcommand{\N}{\mathbb{N}}
  
  \newcommand{\Z}{\mathbb{Z}}

  \newcommand{\cF}{\mathcal{F}}

  \newcommand{\rE}{\mathrm{E}}

  \newcommand{\hs}{\hspace*{\parindent}}
  \newcommand{\proof}{\hs \textbf{Proof.\ }}

  \newcommand{\qed}{\hspace*{\fill} $\Box$\\}

  \newtheorem{theo}{\bfseries \hs Theorem}[section]
  \newtheorem{defn}[theo]{\bfseries \hs Definition}
  \newtheorem{prop}[theo]{\bfseries \hs Proposition}
  \newtheorem{lemma}[theo]{\bfseries \hs Lemma}
  \newtheorem{corol}[theo]{\bfseries \hs Corollary}

  \numberwithin{equation}{section} 

\begin{document}

 \title{Exact conditions for countable inclusion-exclusion identity and extensions}

 \author{
  Shmuel Friedland\\
  \texttt{friedlan@uic.edu}
                   \and
  Elliot Krop
  \\
  \texttt{ekrop1@math.uic.edu}
  }

  \date{Department of Mathematics, Statistics, and Computer Science,\\
        University of Illinois at Chicago\\
        Chicago, Illinois 60607-7045, USA}

 \maketitle

 \begin{abstract}
     We give simple necessary and sufficient conditions for the inclusion-exclusion
     identity to hold for an infinite countable
     number of sets.  In terms of a random variable, whose range
     are nonnegative integers, this condition is equivalent
     to the convergence to zero of binomial moments.  Some
     standard extensions of the countable inclusion-exclusion
     identity are also given.
     \\[\baselineskip] 2000 Mathematics Subject
     Classification: 05A19, 05A20, 60C05\\[\baselineskip]

\end{abstract}

 \section{Introduction}

 The method of inclusion and exclusion and the accompanying
 Bonferroni-type inequalities are very useful and versatile tools in
 probability and combinatorics.  See for example the relative recent book of
 Galambos and Simonelli \cite{GS} and references therein.

 Let $(\Omega,\Sigma,\Pr)$ be a probability space consisting of
 the sample space $\Omega$, a $\sigma$-algebra $\Sigma$ and a
 probability measure $\Pr$ on $\Sigma$.  Let $A_1,\ldots,A_n\in \Sigma$.
 Then the \emph{Inclusion-Exclusion Identity} states:
 \begin{equation}\label{inexp1}
 \Pr(\cup_{i=1}^n A_i)=\sum_{k=1}^n
 (-1)^{k-1}\sum_{1\le i_1<\ldots<i_k\le n}
 \Pr(A_{i_1}\cap\ldots\cap A_{i_k}).
 \end{equation}
 The main purpose of this paper is to give simple necessary and
 sufficient conditions for the \emph{Countable Inclusion-Exclusion Identity}
 \begin{equation}\label{inexc}
 \Pr(\cup_{i\in\N} A_i)=\sum_{k\in\N}
 (-1)^{k-1}\sum_{1\le i_1<\ldots<i_k}
 \Pr(A_{i_1}\cap\ldots\cap A_{i_k}),\;\emph{where } A_i\in\Sigma \textrm{ for
 }i\in\N.
 \end{equation}
 (Here $\N,\Z_+$ is the set of positive and nonnegative integers respectively.)
 As in \cite{GS}, let
 \begin{equation}\label{skdef}
 S_k:=\sum_{1\le i_1<\ldots<i_k}
 \Pr(A_{i_1}\cap\ldots\cap A_{i_k})\in [0,\infty] \textrm{ for }
 k\in\N.
 \end{equation}
 Clearly, for (\ref{inexc}) to hold we must assume that each $S_k$
 is finite, and the sequence $S_k,k\in\N$ converges to $0$.  The
 result of Takacs \cite{Tak} claims that (\ref{inexc}) holds if
 $S_k,k\in\N$ is a sequence of nonnegative numbers that converge
 exponentially to zero:
 \begin{equation}\label{tak}
 S_k\in [0,\infty) \textrm{ for } k\in\N \textrm{ and }\limsup_{k\to\infty} S_k^{\frac{1}{k}}<1.
 \end{equation}
 Moreover, it follows from \cite[p. 111, (39)]{Tak}, that
 the above conditions yield the following generalization
 of (\ref{inexc})
 \begin{equation}\label{inexc1}
 \Pr(\cup_{1 \le i_1<\ldots <i_k}(A_{i_1}\cap\ldots\cap
 A_{i_k}))=\sum_{j\in\Z_+} (-1)^{j}{j+k-1\choose k-1}S_{j+k}
 \textrm{ for all } k\in\N.
 \end{equation}
 In this note we show
 \begin{theo}\label{inexmt} Let $(\Omega,\Sigma,\Pr)$ be
 probability space and assume that $A_i\in\Omega$ for $i\in\N$.
 Then (\ref{inexc1}) holds for some $k\in\N$ if and only if
 $S_l,l\in\N$ is a sequence of nonnegative numbers such that
 $\lim_{l\to\infty} l^{k-1} S_l=0$.
 \end{theo}
 We now list briefly the contents of this paper.
 In \S2 we prove Theorem \ref{inexmt}.  In \S3 we discuss a
 partition of $\Omega$ induced by $A_i,i\in\N$ and discuss
 a few applications.  In \S4 we give an analog of
 Theorem \ref{inexmt} for the random variable $X:\Omega \to \Z_+$.

 \section{Proof of the main theorem}
 Let the assumptions of Theorem \ref{inexmt} hold.
 Denote
 \begin{equation}\label{skndef}
 S_{k,n}:=\sum_{1\le i_1<\ldots<i_k\le n}
 \Pr(A_{i_1}\cap\ldots\cap A_{i_k}) \textrm{ for }
 k\le n \textrm{ and } S_{k,n}:=0 \textrm{ for } k >n.
 \end{equation}
 Clearly, $0\le S_{k,n}\le {n\choose k}$ for each $k,n\in\N$.
 \begin{lemma}\label{limskn}  Let the assumptions of Theorem \ref{inexmt}
 hold.  Then for each $k\in\N$ the sequence $S_{k,n},n\in\N$ is a
 nondecreasing sequence that converges in the generalized sense to
 $S_k\in[0,\infty]$.  That is, $S_k<\infty$ if and only
 $S_{k,n},n\in\N$ is a bounded sequence converging to $S_k$, and
 $S_k=\infty$ if and only if $S_{k,n},n\in\N$ is an unbounded
 sequence.
 \end{lemma}
 The proof of this lemma is standard and is left to the reader.

 \textbf{Proof of Theorem} \ref{inexmt}.
 Fix $k\in\N$.  Assume first that (\ref{inexc1}) holds.
 First, $S_j,j\in\N$ is a sequence of nonnegative numbers.
 Second, the convergence the series in (\ref{inexc1}) implies
 that $\lim_{l\to\infty} {l-1\choose k-1} S_l=0$, which is
 equivalent to $\lim_{j\to\infty} l^{k-1} S_l=0$.

 Assume first that $S_l,l\in\N$ is a sequence of nonnegative
 numbers.
 Recall the Bonferroni inequalities \cite[Ineq. I.2]{GS}.
 Let $1\le k\le n$ and assume that $d,r\in\Z_+$.
 Then
 \begin{eqnarray}
 \sum_{j=0}^{2d+1} (-1)^{j}{j+k-1\choose k-1}S_{j+k,n}\le
 \Pr(\cup_{i_1<\ldots <i_k\le n}(A_{i_1}\cap\ldots\cap
 A_{i_k}))\le \nonumber\\
 \sum_{j=0}^{2r} (-1)^{j}{j+k-1\choose k-1}S_{j+k,n}.
 \label{bonfer1}
 \end{eqnarray}
 Let $n\to\infty$.  Clearly, $\Pr(\cup_{i_1<\ldots <i_k\le n}(A_{i_1}\cap\ldots\cap
 A_{i_k}))\nearrow \Pr(\cup_{i_1<\ldots <i_k}(A_{i_1}\cap\ldots\cap
 A_{i_k}))$.  Use Lemma \ref{limskn} to deduce the Bonferroni type
 inequality
 \begin{eqnarray}
 \sum_{j=0}^{2d+1} (-1)^{j}{j+k-1\choose k-1}S_{j+k}\le
 \Pr(\cup_{i_1<\ldots <i_k}(A_{i_1}\cap\ldots\cap
 A_{i_k}))\le \nonumber\\
 \sum_{j=0}^{2r} (-1)^{j}{j+k-1\choose k-1}S_{j+k}.
 \label{bonferc1}
 \end{eqnarray}
 Let $a_d \le c_k\le b_r$ be the left-hand side, the middle part and the right-hand
 side of the above inequalities.  Then $c_k\in [a_d,b_d]$ for any
 $d\in\Z_+$.  Assume second that $\lim_{l\to\infty} {l\choose k-1} S_l=0$.
 Hence $b_d-a_d={2d+k \choose k-1} S_{2d+k+1}\to 0$.  Therefore
 the left-hand side and the right-hand side of (\ref{bonferc1})
 converge to $c_k$.  \qed

 \begin{corol}\label{corthm1}  Let $(\Omega,\Sigma,\Pr)$ be
 a probability space, assume that $S_l\in [0,\infty),l\in\N$ and
 (\ref{inexc1}) holds for some $k=m >1$.  Then (\ref{inexc1}) holds for
 $k=1,\ldots,m-1$.
 \end{corol}

 \section{A decomposition of a countable sets}
 \begin{defn}\label{defdecom}  Let $\Omega$ be an infinite set and
 $A_i\subseteq \Omega,i\in\N$.  Then
 \begin{itemize}
 \item Let $\cF\subset 2^{\Omega}$ be the set of all
 nonempty finite subsets of $\Omega$, let
 $\tilde\cF:=\cF\cup\{\emptyset\}\subset 2^{\Omega}$ and for each $j\in\Z_+$
 let $\cF_j\subset \tilde\cF$ be the set all finite subsets of $\N$ of cardinality $j$.
 \item $B_{\emptyset}:=\Omega\backslash \cup_{i\in\N} A_i$.
 \item $B_{\infty}:=\lim\sup A_i=\cap_{i\in\N}\cup_{j\ge i} A_j$ the set of
 points that belong to an infinite number of $A_i,i\in\N$.
 \item For each $U\in \cF$ denote by
 $A_U:=\cap_{i\in U} A_i$, and by $B_U:=A_U\backslash \cup_{i\in \N\backslash U}
 A_{U\cup\{i\}}$ the set of points belonging only to $A_i,i\in U$.
 Let $A_U':=A_U\backslash B_{\infty}$.

 \end{itemize}
 \end{defn}

 \begin{prop}\label{propset}  Let $\Omega$ be an infinite set and
 $A_i\subseteq \Omega,i\in\N$.  Then
 \begin{enumerate}
 \item $B_{\infty}$ and
 $B_U,U\in\tilde\cF$ form a countable partition
 of $\Omega$.
 \item
 $\cup_{i\in\N} A_i=B_{\infty}\cup(\cup_{U\in\cF}
 B_U)$.
 \item $A_U=A'_U \cup (A_U\cap B_{\infty})$.

 \end{enumerate}
 Assume in addition that $(\Omega,\Sigma,\Pr)$ is a probability
 space and $A_i\in \Sigma, i\in \N$.  If $S_k$ is finite for some
 $k\in\N$ then $\Pr(B_{\infty})=0$.  Let
 \begin{equation}\label{deftk}
 T_j:=\sum_{U\in\cF_j} \Pr (B_U), \quad \textrm{for each
 }j\in\Z_+.
 \end{equation}
 Then $\sum_{j\in\Z_+} T_j\le 1$ and equality holds if and only if $\Pr(B_{\infty})=0$.
 Assume that $\Pr(B_{\infty})=0$.  Then
 \begin{equation}\label{sktkid}
 S_k=\sum_{j\in\Z_+} {j+k \choose k} T_{j+k} \in [0,\infty], \quad \textrm{for each
 }k\in\N.
 \end{equation}
 Let $l\in\N$ and assume that $S_l\in [0,\infty)$.  Then for each
 positive integer $k< l$ $S_k\in [0,\infty)$.  Suppose furthermore
 that $2d+k+1,2r+k\in
 [1,l]$ for some $d,r\in\Z_+,k\in\N$.  Then (\ref{bonferc1})
 holds.
 \end{prop}

 \proof  Claims \textit{1, 2, 3} of the proposition are
 straightforward.  Assume that $(\Omega,\Sigma,\Pr)$ is a probability
 space and $A_i\in \Sigma, i\in \N$.  Since $\cF$ is countable, $A_U,A'_U,B_{\infty}\in \Sigma$
 for each $U\in \tilde\cF$.

 Assume that $S_k=\sum_{U\in\cF_k} \Pr(A_U)<\infty$.
 Let $B_{k,\infty}$ be the set of elements in $\Omega$
 which belong to an infinite number of $A_U$, where $U\in\cF_k$.
 Use the Borel-Cantelli Lemma to deduce that
 $\Pr(B_{k,\infty})=0$.
 It is straightforward to show that $B_{k,\infty}=B_{\infty}$.  Hence
 $\Pr(B_{\infty})=0$.

 From \textit{1} it follows that $\Pr(B_{\infty})+\sum_{j\in\Z_+}
 T_j=1$.  Hence  $\sum_{j\in\Z_+} T_j\le 1$ and equality holds if and only if $\Pr(B_{\infty})=0$.

 Assume that $\Pr(B_{\infty})=0$.  Then $\Pr(A_U)=Pr(A_U')$ for
 $U\in\cF$ and $S_k=\sum_{U\in\cF_k} \Pr(A_U')$.
 Let $V\in \cF_l$, where $l\ge k$.  Then $B_V\subset A_U'$ for each
 $U\subset V$.  Suppose that $\#U=k$.  Then $V$ has exactly
 $l\choose k$ distinct $k$-elements subsets $U$.  Thus $B_V$ is a
 subset of exactly $l\choose k$ sets $A_U', U\in \cF_k$, and all
 other subsets $A_U',U\in\cF_k$ are disjoint from $B_V$.  Hence
 (\ref{sktkid}) holds in the generalized sense, i.e. $S_k=\infty$ if and
 only if the right-hand side of (\ref{sktkid}) diverges.

 Suppose that $S_l<\infty$.  Then $\Pr(B_{\infty})=0$.
 Furthermore the series (\ref{sktkid}) converges for $k=l$.
 Let $l > k\in\N$.  As ${p\choose l}> {p\choose k}$ for $p \ge
 2l-1$ it follows that
 $$S_l\ge\sum_{p\ge 2l-1} {p \choose l} T_l \ge \sum_{p\ge 2l-1} {p \choose k}
 T_l=S_k -\sum_{k\le p\le 2l-2} {p \choose l} T_l.$$
 Hence $S_k<\infty$.

 Suppose furthermore that $2d+k+1,2r+k\in
 [1,l]$ for some $d,r\in\Z_+,k\in\N$.  Then the proof of Theorem
 \ref{inexmt} yields (\ref{bonferc1}).  \qed

 \section{Random variables with values in $\Z_+$}
 As in \cite{GS}, Theorem \ref{inexmt} or its variation can be
 reformulate in terms of binomial moments of a random variable
 $X:\Omega \to \Z_+$.  Assume that $\Pr(X=j)=T_j$ for each $j\in
 \Z_+$.  Hence $\sum_{j\in\Z_+} T_j=1$.  Let $S_j:=\rE({X \choose
 j})$ for $j\in \Z_+$.  Then $S_0=1$ and $S_k$ is given by
 (\ref{sktkid}) for $k\in\N$.

 \begin{theo}\label{binmom}  Let $(\Omega,\Sigma,\Pr)$ be
 probability space and assume that $X:\Omega\to\Z_+$ is a random
 variable.  Let $S_j:=\rE({X \choose
 j})\in [0,\infty]$ for $j\in \Z_+$.  Suppose that $S_l<\infty$.
 Then $S_k<\infty$ for each $l > k\in\N$.
 Suppose furthermore
 that $2d+k+1,2r+k\in
 [1,l]$ for some $d,r\in\Z_+,k\in\N$.
 Then
 \begin{eqnarray}
 \sum_{j=0}^{2d+1} (-1)^{j}{j+k-1\choose k-1}S_{j+k}\le
 \Pr(X\ge k) \le
 \sum_{j=0}^{2r} (-1)^{j}{j+k-1\choose k-1}S_{j+k},
 \label{bonferc2}\\
 \sum_{j=0}^{2d+1} (-1)^{j}{j+k-1\choose k-1}S_{j+k-1}\le
 \Pr(X= k-1) \le
 \sum_{j=0}^{2r} (-1)^{j}{j+k-1\choose k-1}S_{j+k-1}.
 \label{bonferc3}
 \end{eqnarray}
 Furthermore for $k\in\N$ one has the equalities
 \begin{eqnarray}
 \Pr(X\ge k) =
 \sum_{j\in\Z_+} (-1)^{j}{j+k-1\choose k-1}S_{j+k},
 \label{binidc2}\\
 \Pr(X= k-1) =
 \sum_{j\in\Z_+} (-1)^{j}{j+k-1\choose k-1}S_{j+k-1},
 \label{binidc3}
 \end{eqnarray}
 if and only if $S_l,l\in\N$ is a sequence of nonnegative numbers such that
 $\lim_{l\to\infty} l^{k-1} S_l=0$.
 \end{theo}
 \proof  Assume that $n_0\in\Z_+$ is the first nonnegative integer
 such that $\Pr(X\le n_0) >0$.  For $\N\ni n\ge n_0$ let $X_n:\Omega\to
 [0,n]\cap \Z_+$ be the random variable whose distribution is
 given by $\Pr(X_n=k)=\frac{\Pr(X=k)}{\Pr (X\le n)}$ for
 $k=0,\ldots,n$.  Let $S_{j,n}:=\rE({X_n\choose j})$ for $\in
 \Z_+$.  Then
 $$\lim_{n\to\infty} S_{j,n}=S_j,\;\lim_{n\to\infty} \Pr(X_n\ge
 j)=\Pr(X\ge j),\; \lim_{n\to\infty} \Pr(X_n=j)=\Pr(X= j),\;\textrm{for all
 }j\in\N.$$
 Use the two type of Bonferroni inequalities given in
 \cite[Ineq. I.2]{GS} and the arguments of the proof of Theorem
 \ref{inexmt} and Proposition \ref{propset} to deduce the theorem.
 \qed

 If $S_k,k\in\N$ is a sequence of nonnegative numbers such that
 $\limsup S_k^{\frac{1}{k}}<1$ then the equalities
 (\ref{binidc2}-\ref{binidc3}) are due to Takacs \cite{Tak}.


\begin{thebibliography}{99}
 \bibitem{GS} J. Galambos and I. Simonelli, \emph{Bonferroni-type Inequalities with
 Applications}, Springer, 1996.
 \bibitem{Tak} L. Takacs, On the Method of Inclusion and
 Exclusion, \emph{J. Amer. Stat. Assn.} 62 (1967), 102-113.
 \end{thebibliography}
\end{document}